\documentclass[11pt]{amsart}

\usepackage{amssymb}
\usepackage{amsfonts}
\usepackage{amsmath}
\usepackage{graphicx}
\usepackage{amscd}
\usepackage{enumerate} 
\newtheorem{theorem}{Theorem}
\numberwithin{theorem}{section}

\newtheorem{case}{Case}

\newtheorem{corollary}[theorem]{Corollary}

\newtheorem{lemma}[theorem]{Lemma}

\newtheorem{proposition}[theorem]{Proposition}
\newtheorem{remark}[theorem]{Remark}

\theoremstyle{definition}
\newtheorem{claim}{Claim}

\newenvironment{proof1}[1][Proof]{\noindent\textbf{#1.} }{\qed}

\DeclareMathOperator{\dist}{dist}
\DeclareMathOperator{\lo}{long}
\DeclareMathOperator{\Tr}{Trans}
\DeclareMathOperator{\id}{id}

\DeclareMathOperator{\diam}{diam}

\DeclareMathOperator{\ord}{ord}
\DeclareMathOperator{\cl}{cl}
\DeclareMathOperator{\Int}{int}
\DeclareMathOperator{\ext}{ext}
\DeclareMathOperator{\bd}{bd} 

\usepackage[normalem]{ulem}

\makeatletter
\@namedef{subjclassname@2020}{%
  \textup{2020} Mathematics Subject Classification}
\makeatother

\let\oldproofname=\proofname
\renewcommand{\proofname}{\rm\bf{\oldproofname}}

\begin{document}

\title{The hyperspace $\omega(f)$ when $f$ is a transitive dendrite mapping}

\author[J. M. Mart\'inez-Montejano]{Jorge M. Mart\'inez-Montejano}
\address{(J. M. Mart\'{\i}nez-Montejano) Departamento de Matem\'{a}ticas,
Facultad de Ciencias, Universidad Nacional Aut\'{o}noma de M\'{e}xico,
Circuito Exterior, Cd. Universitaria,  Ciudad de M\'{e}xico, 04510, M\'{e}xico}
\email{jorgemm@ciencias.unam.mx}

\author[H. M\'endez]{H\'ector M\'endez}
\address{(H. M\'endez) Departamento de Matem\'{a}ticas,
Facultad de Ciencias, Universidad Nacional Aut\'{o}noma de M\'{e}xico,
Circuito Exterior, Cd. Universitaria,  Ciudad de M\'{e}xico, 04510, M\'{e}xico}
\email{hml@ciencias.unam.mx}

\thanks{This paper was partially supported by the project ``Sistemas din\'amicos discretos y teor\'{\i}a de continuos I" (IN105624) of PAPIIT, DGAPA, UNAM}

\author[Y. N. Vel\'azquez-Inzunza]{Yajaida N. Vel\'azquez-Inzunza}
\address{(Y. N. Vel\'azquez-Inzunza) Departamento de Matem\'{a}ticas,
Facultad de Ciencias, Universidad Nacional Aut\'{o}noma de M\'{e}xico,
Circuito Exterior, Cd. Universitaria,  Ciudad de M\'{e}xico, 04510, M\'{e}xico}
\email{ynvi@ciencias.unam.mx}


\begin{abstract}
Let $X$ be a compact metric space. By $2^X$ we denote the hyperspace of all closed and non-empty subsets of $X$ endowed with the Hausdorff metric. Let $f:X\to X$ be a continuous function. In this paper we study some topological properties of the hyperspace $\omega(f)$, the collection of all omega limits sets $\omega(x,f)$ with $x\in X$. We prove the following: $i)$ If $X$ has no isolated points, then, for every continuous function $f:X\to X$, $\Int_{2^X}(\omega(f))=\emptyset$. $ii)$ If $X$ is a dendrite for which every arc contains a free arc and $f:X\to X$ is transitive, then the hyperspace $\omega(f)$ is totally disconnected. $iii)$ Let $D_\infty$ be the Wazewski's universal dendrite. Then  there exists a transitive continuous function $f:D_\infty\to D_\infty$ for which the hyperspace  $\omega(f)$ contains an  arc; hence, $\omega(f)$ is not totally disconnected.

\end{abstract}

\date{\today}
\maketitle

\section{Introduction }

A \textit{map} is a continuous  function. Let $\mathbb{N}$
denote the set of positive integers. 
Let $X$ be a compact metric space. Given $p \in X$ and $\varepsilon > 0$, let $B(p,\varepsilon)$ denote
the $\varepsilon$-ball in  $X$ around $p$. Also, for a given subset $A$ of
$X$, the symbols $\cl_{X}(A)$, $\Int_{X}(A)$ and $\bd_{X}(A)$ denote the
closure, the interior and the boundary of   the set $A$ in $X$, respectively.

A compact metric space $X$ is said to be \textit{totally disconnected} if
for each pair of distinct points $x,y \in X$, there exist open subsets $U$
and $V$ of $X$ such that $x \in U$, $y \in V$, $U\cap V=\emptyset$ and $X=U\cup
V$.

A \textit{dynamical system} is a pair $(X,f)$ where $X$ is a compact
metric space and $f:X\to X$ is a map.
Given a dynamical system $(X,f)$, we denote by $f^{n}$ the $n$-th iterate
of $f$; that is, $f^{0} = \id_{X}$ (the identity map) and $f^{n}= f \circ
f^{n-1}$ if $n \geq 1$. For each $x\in X$, the \emph{orbit of $x$ under $f$} is the sequence $o(x,f)=\lbrace x,f(x),f^{2}(x),\ldots \rbrace$.
A point $x \in X$ is called \textit{periodic} if there is \(k\in\mathbb{N}\)
such that $f^{k}(x)=x$. Given a periodic point $x\in X$, the \textit{period}
of $x$ is the least positive integer $m$ such that
$f^{m}(x)=x$; when $m=1$ we say that $x$ is a \textit{fixed point}.
A non-empty subset $A$ of $X$ is called \textit{strongly invariant under f} provided that
$f(A) = A$.
We define the $\omega$-\textit{limit set} of a point $x\in X$ to be the set
$$\omega(x,f)=\bigcap\limits_{n\in \mathbb{N}} \cl_{X}(\{ f^{k}(x): k \geq
n \}).$$
Note that $y\in\omega(x,f)$ if and only if there exists a sequence of natural
numbers  $n_{1} < n_{2}<\ldots$ such that $\lim\limits_{i \to \infty}f^{n_{i}}(x)=y$.
It is well known that $\omega(x,f)$ is a non-empty, closed and strongly invariant
subset of  $X$ (see \cite{Block}).
The map $f:X \rightarrow X$ is 
\begin{itemize}
\item[(i)] \emph{transitive} provided that for each pair of non-empty open subsets $U$ and $V$ of $X$, there exists $n\in \mathbb{N}$ such that $f^{n}(U)\cap V\neq \emptyset$. This condition is equivalent to the existence of a point $x\in X$ with $\omega(x,f)=X$ (see \cite{Block});
\item[(ii)] \emph{mixing} provided that for each pair of non-empty open sets $U$ and $V$ of $X$, there exists $N\in \mathbb{N}$ such that for every $n\geq N$, $f^{n}(U)\cap V\neq \emptyset$.
\end{itemize}

A \textit{continuum} is a non-empty compact connected metric space with more than one point. A \textit{subcontinuum} of a continuum $X$ is a non-empty closed connected subset of $X$. 
An \textit{arc} is a space homeomorphic to the unit interval $[0,1]$. A \textit{graph} is a continuum that can be expressed as the finite union of arcs where each pair of them meet in a subset of their end points.

A \textit{dendrite} is a locally connected continuum containing no simple
closed curves. 
For a dendrite $D$ it is known that every connected subset of $D$ is arcwise
connected  (\cite[Proposition 10.9, p. 169]{Cont}), and that the intersection
of any two connected subsets of $D$ is connected (\cite[Theorem 10.10, p.169]{Cont}).

If $X$ is a continuum and $p\in X$, then the (Menger-Urysohn) \textit{order}
of $p$ in $X$, denoted by \(\ord_X(p)\), is defined to be the least cardinal
\(\kappa\) such that $p$ has arbitrarily small neighborhoods with boundary
of cardinality less or equal than \(\kappa\). For a dendrite $D$ and $p\in
D$, \(\ord_X(p)\) is equal to the number of components of $X\setminus\{p\}$.
A point in a dendrite is an \textit{end point} if it has order 1, an \textit{ordinary
point} if it has order 2 and a \textit{ramification point} if it has order
at least 3. Given a dendrite $D$, the set of end points, ordinary points and
ramification points are denoted by $E(X)$, $O(X)$, and $R(X)$, respectively. 

Let $D$ be a dendrite. Let $a,b \in D$, with $a\neq b$. With $[a,b]$ we denote the unique arc contained in $D$ that has $a$ and $b$ as end points. The symbol $(a,b)$ denotes the set $[a,b]\setminus \lbrace a, b \rbrace$.
Notice that, for each $c\in(a,b)$, $c\notin E(X)$. An arc $[p,q]$ in $D$ is called a \textit{free arc} provided that $(p,q)$ is an open subset of $D$ or, equivalently, if $(p,q)$ contains no ramification points of $D$.

Given a compact metric space $X$, with metric $d$, we define the \textit{hyperspace of compact subsets of $X$} as
\[2^{X}= \{A\subseteq X: A\textnormal{ is closed and non-empty}\}.\]
This hyperspace is endowed with the Hausdorff metric $H$ (see \cite[Definition 4.1, p.52]{Cont}). 

Given a finite collection $A_{1},\ldots,A_{n}$ of subsets of $X$, define
\[\langle A_{1},\ldots,A_{n}\rangle=\{K\in2^{X}:K \subseteq A_{1}\cup\cdots \cup A_{n}\textnormal{ and, for each }i\in\{1,\ldots,n\}\text{, }  K\cap A_{i}\neq\emptyset\}.\]
It is known that the family of subsets of $2^{X}$ of the form $\langle U_{1}, \ldots, U_{n} \rangle$, where $U_{1},\ldots,U_{n}$ is a finite collection of open subsets of $X$, is a basis for the topology generated by the Hausdorff metric (see \cite[Theorem 4.5, p. 54 ]{Cont}). Also, it is easy to see that if $F_{1},\ldots,F_{n}$ is a finite collection of closed subsets of $X$, then $\langle F_{1},\ldots, F_{n} \rangle$ is a closed subset of $2^X$.

Given a dynamical system $(X,f)$, define \textit{the hyperspace of $\omega$-limit sets of $f$} as $$\omega(f) = \{ A\in 2^X:\text{there exists }x\in X \text{ with }A=\omega(x,f)\}.$$

The topological properties of $\omega(f)$  (as a subspace of $2^X$) has been widely studied; for instance, even thought in general it is not true that $\omega(f)$ is compact \cite[Example 2.2]{Mai}, it is known that if $(G,f)$ is a dynamical system where $G$ is a graph, then $\omega(f)$ is compact \cite[Theorem 3.]{Mai}. For dendrites, there are examples of maps for which its hyperspace of $\omega$-limit sets is not compact (see, e.g., \cite[Theorem 2]{Kocan}) but it is known that if $(D,f)$ is a dynamical system where $D$ is a dendrite and $f$ is a monotone map (that is, for each $x\in D$, $f^{-1}(x)$ is
 connected), then $\omega(f)$ is compact \cite[Theorem 4.5]{abouda2016}.

In this paper we contribute to this study. First, in section 3, we prove that if $(X,f)$ is a dynamical system where $X$ has no isolated points, then $\Int_{2^X}(\omega(f))=\emptyset$. Then, we show that $\omega(f)$ is totally disconnected whenever $(D,f)$ is a dynamical system where $D$ is a dendrite for which every arc contains a free arc and $f$ is transitive. Finally, in section 4, we describe a mixing map $f$ on the  Wazewski's universal dendrite for which $\omega(f)$ is not totally disconnected.

\section{ Previous Results }

Let $(X, f)$ be a dynamical system. The following two results are well known.
\begin{remark}\label{conjunto_de_omegas_denso}\cite[Theorem 2.8]{Akin}
If $f:X\to X$ is transitive, then the set $\Tr(f)=\{x\in X:\cl_X(o(x,f))=X\}$ is residual. Hence, $\Tr(f)$ is a dense subset of $X$. For each $x\in\Tr(f)$ we have that $\omega(x,f)=X$. Thus, the set $\{x\in X:\omega(x,f)=X\}$ is a dense subset of $X$. 
\end{remark}

\begin{lemma}\label{omega finito}\cite[Lemma 4]{Block}
Let  $x\in X$. Then  $\omega(x,f)$
contains only finitely many points if and only if there exists a periodic
point $z\in X$ such that $\omega(x,f)=o(z,f)$.
\end{lemma}

Also, it is well known that dendrites are locally connected regular curves, we express this fact in the following lemma for future reference.

\begin{lemma}\label{lema_base_conectos_fr_finita}\cite[Theorem 2.20]{Cont}
Let $D$ be a dendrite. Then $D$ has a basis $\mathcal{B}$ such that each
member of $\mathcal{B}$ is a connected open subset of \(X\) with finite boundary.
\end{lemma}

The following lemma is not difficult to prove.

\begin{lemma} \label{cerradura_vietorico}
Let $X$ be a compact metric space and let $U_{1},\ldots,U_{n}$ be a finite collection of subsets of $X$. Then $\cl_{2^{X}}( \langle U_{1}, \ldots, U_{n} \rangle ) =\langle \cl_{X}(U_{1}), \ldots, \cl_{X}(U_{n}) \rangle $.
\end{lemma}

\begin{lemma} \label{elemento_en_la_frontera}
Let $X$ be a compact metric space and let $U_{1},\ldots, U_{n}$ be a finite collection of  non-empty subsets of $X$. Consider $\mathcal{U} =\langle U_{1},\ldots,U_{n}\rangle $ and let $A\in 2^X$. If $A \in \cl_{2^{X}}(\mathcal{U})\setminus\mathcal{U}$, then  there exist $a\in A$ and $j\in\{1,\ldots,n\}$ such that $a\in\cl_{X}(U_{j}) \setminus U_{j}$. 
\end{lemma}

\begin{proof1}
By Lemma \ref{cerradura_vietorico}, $A\in\langle \cl_{X}(U_{1}),\ldots,\cl_{X}(U_{n}) \rangle\setminus\langle U_{1}, \ldots, U_{n} \rangle$. We consider two cases.
        
\begin{case}
There exists $j\in \{1,\ldots,n\}$ such that $A\cap U_{j} =\emptyset$.
\end{case} 
Since $A \cap \cl_{X}(U_{j})\neq\emptyset$, there exists $a \in A\cap\cl_{X}(U_{j})$. Hence, $a\in\cl_{X}(U_{j}) \setminus U_{j}$.

\begin{case}                
$A\nsubseteq\bigcup\limits_{i=1}^{n} U_{i}$.
\end{case}
Then, there exists $a \in A$ such that $a \notin \bigcup \limits_{i=1}^{n} U_{i}$. Since $A\subseteq\bigcup\limits_{i=1}^{n} \cl_X(U_{i})$, there is $j\in\{1,\ldots,n\}$ such that $a\in\cl_{X}(U_{j})$. We have that $a\in\cl_{X}(U_j)\setminus U_{j}$.
\end{proof1}


\section{Results}

\begin{theorem}
Let $(X,f)$ be a dynamical system. If $X$ has no isolated points, then $\Int_{2^{X}}(\omega(f)) = \emptyset$.
\end{theorem}

\begin{proof1}
Assume, on the contrary, that $\Int_{2^{X}}(\omega(f)) \neq \emptyset$. Take $A\in\Int_{2^{X}}(\omega(f))$ and $U_1,\dots,U_k$ a finite collection of non-empty open subsets of $X$ such that $A\in\langle U_1,\ldots,U_k\rangle\subseteq\omega(f)$. Choose $x_1,\ldots,x_k$ pairwise distinct points in $X$ such that, for each $j\in\{1,\dots,k\}$, $x_j\in U_j$. Let $B=\{x_1,\ldots,x_k\}$. We have that $B\in\langle U_1,\ldots,U_k\rangle\subseteq\omega(f)$. Since $B$ is a finite subset of $X$, by Lemma \ref{omega finito}, there is a periodic point $z\in X$ such that $B=o(z,f)$. Let $y\in X$ be such that $y\in U_1$ and  $y\notin B$. Define $C=B\cup\{y\}$. Note that $C\in\langle U_1,\ldots,U_k\rangle\subseteq\omega(f)$. Since $C$ is a finite subset of $X$, by Lemma \ref{omega finito}, there is  a periodic point $u\in X$ such that $C=o(u,f)$. We have that $z\in B=o(f,z)$, where $z$ is a periodic  point of $f$ and the period of $z$ is $k$; also, $z\in C=o(u,f)$, where $u$ is a periodic point of $f$ and the period of $u$ is $k+1$, which is a contradiction.

Therefore, $\Int_{2^{X}}(\omega(f)) = \emptyset$.
\end{proof1}     

\begin{theorem} \label{teorema_w(f)_totalmete_disconexo}
Let $(D,f)$ be a dynamical system where $D$ is a dendrite. If for each arc $[a,b]$ contained in $D$, with $a\neq b$, there exists $c\in (a,b)$ such that $\omega(c,f)=D$, then $\omega(f)$ is a totally disconnected subspace of $2^{D}$.
\end{theorem}

\begin{proof1}
Let $\omega(x,f)$ and $\omega(y,f) $ be two distinct elements of $\omega(f)$. We may assume, without loss of generality, that there exists $p\in\omega(x,f) \setminus\omega(y,f)$. Let $\delta >0$ be such that $B(p,\delta)\cap\omega(y,f)=\emptyset$. By Lemma \ref{lema_base_conectos_fr_finita}, there is an open connected subset $V$ of $X$ such that $p\in V\subseteq\cl_D(V)\subseteq B(p,\delta)$ and $\bd_D(V)=\{a_1,\dots,a_k\}$. Consider, for each $i\in\{1,\ldots,k\}$, the arc $[p,a_i]$. By hypothesis, for each $i\in\{1,\ldots,k\}$, there is a point $r_i\in (p,a_i)$ such that $\omega(r_i,f)=D$. 

Given $i\in\{1,\ldots,k\}$,  since $r_i\notin E(X)$, we have that $D\setminus\{r_i\}$ is disconnected \cite[Theorem 10.7]{Cont}; hence, we can write $D\setminus\{r_i\}=H_i\cup K_i$ where $H_i$ is the component of $D\setminus\{r_i\}$ that contains $p$. Note that $a_i\notin H_i$; thus, since $a_i\neq r_i$, $a_i\in K_i$. 

For each $i\in\{1,\ldots,k\}$, $H_i$ and $K_i$ are open subsets of $D$, and $\cl_{D}(H_i)=H_i\cup\{r_i\}$ and $\cl_{D}(K_i)=K_i\cup\{r_i\}$ are proper subcontinua of $D$ \cite[Proposition 6.3]{Cont}.

We have that $\bigcap\limits_{i=1}^k\cl_{D}(H_i)$ is a connected subset of $D$ \cite[Theorem
10.10]{Cont} such that $$p\in\bigcap\limits_{i=1}^k\cl_D(H_i)\subseteq D\setminus\{a_1,\ldots,a_k\}=
D\setminus\bd_D(V)=\Int_D(V)\cup\ext_D(V)=V\cup(D\setminus\cl_{D}(V)).$$ Since $p\in V$, $\bigcap\limits_{i=1}^k\cl_D(H_i)\subseteq
V$.

We show that there must be a $j\in\{1,\ldots,k\}$ such that $\omega(y,f)\cap K_j\neq\emptyset$. Assume, on the contrary, that for each $i\in\{1,\ldots,k\}$, $\omega(y,f)\cap K_i=\emptyset$. Then $\omega(y,f)\subseteq\bigcap\limits_{i=1}^k\cl_D(H_i)\subseteq V$; a contradiction.

Now, define $B=\{j\in\{1,\ldots,k\}:\omega(y,f)\cap K_j\neq\emptyset\}$. Let us denote $B=\{j_1,\ldots,j_m\}$ and consider $$\mathcal{U}=\langle K_{j_1},\ldots,K_{j_m}\rangle \text{ and }  \mathcal{V}=2^D\setminus\cl_{2^D}(\mathcal{U}).$$ We have that $\mathcal{U}$ and $\mathcal{V}$ are open subsets of $2^D$ such that $\mathcal{U}\cap\mathcal{V}=\emptyset$. We prove the following claims.
\begin{claim}
$\omega(y,f)\in\mathcal{U}$.
\end{claim}
By definition, for each $i\in\{1,\ldots,m\}$, $w(y,f)\cap K_{j_i}\neq\emptyset$. Since $\omega(y,f)\cap V=\emptyset$ and $\bigcap\limits_{i=1}^k\cl_D(H_i)\subseteq
V$, we have that $$\omega(y,f)\subseteq (D\setminus V)\subseteq \left(D\setminus\bigcap\limits_{i=1}^k\cl_{D}(H_i)\right)=\bigcup\limits_{i=1}^k K_i.$$ Thus, $\omega(y,f)\subseteq\bigcup\limits_{i=1}^m
K_{j_i}$. This finishes the proof of Claim 1.

\begin{claim}
$\omega(x,f)\in\mathcal{V}$.
\end{claim}

By Lemma \ref{cerradura_vietorico}, $\cl_{2^D}(\mathcal{U})=\langle\cl_D(K_{j_1}),\ldots,\cl_D(K_{j_m})\rangle$.
 Since $p\in\bigcap\limits_{i=1}^k H_i$, it follows that $p\notin\bigcup\limits_{i=1}^k\cl_D(K_i)$. Thus, since $p\in\omega(x,f)$, we have that $w(x,f)\nsubseteq\bigcup\limits_{i=1}^m\cl_D (K_{j_i})$. Hence, $\omega(x,f)\notin\cl_{2^D}(\mathcal{U})$. This concludes the proof of Claim 2.

\begin{claim}
$\omega(f)\subseteq\mathcal{U}\cup\mathcal{V}$.
\end{claim}

First, note that $p\notin\bigcup\limits_{i=1}^k \cl(K_i)$. Now, assume, on the contrary, that there exist $q\in D$ such that $\omega(q,f)\notin\mathcal{U}\cup\mathcal{V}$. Then  $\omega(q,f)\notin\mathcal{U}$  and $\omega(q,f)\notin\mathcal{V}$; thus, $\omega(q,f)\in\cl_{2^{D}}(\mathcal{U})$. That is, $\omega(q,f)\in\cl_{2^D}(\mathcal{U})\setminus\mathcal{U}$. By Lemma \ref{elemento_en_la_frontera}, there are $e\in\omega(q,f)$ and $s\in\{1,\ldots,m\}$ such that $e\in\cl_D(K_{j_s})\setminus K_{j_s}$. Since $\cl_D(K_{j_s})\setminus K_{j_s}=\{r_{j_s}\}$, we have that $r_{j_s}=e\in\omega(q,f)$. Hence, since $\omega(q,f)$ is closed and strongly invariant, we obtain that $$D=\omega(r_{j_s},f)\subseteq\omega(q,f)\subseteq D.$$  That is, $\omega(q,f)=D$. Since $D=\omega(q,f)\in\cl_{2^D}(\mathcal{U})=\langle\cl_D(K_{j_1}),\ldots,\cl_D(K_{j_m})\rangle$, we have that $D\subseteq\bigcup\limits_{i=1}^m\cl_D
(K_{j_i})$; a contradiction. This finishes the proof of  Claim 3.

Therefore, $\omega(f)$ is totally disconnected.
\end{proof1}     
 
\begin{corollary}\label{corolario_dendrita_totalmente_disconexo}
Let $D$ be a dendrite. Assume that for each arc $[a,b]$ in $D$, with $a\neq b$, there exist $u,v\in [a,b]$, with $u\neq v$, such that $[u,v]$ is a free
arc in $D$.  If $f: D\to D$ is a transitive map, then $\omega(f)$ is totally disconnected.
\end{corollary}

\begin{proof1}
     Let $[a,b]$ be an arc in $D$, with $a\neq b$. By hypothesis, there exist $u,v\in [a,b]$, with $u\neq v$, such that $[u,v]$ is a free arc in $D$. 
Since $f$ is a transitive map, by Remark \ref{conjunto_de_omegas_denso}, there exists a point $c\in(u,v)\subseteq(a,b)$ such that $\omega(c,f) = D$.
Thus, by Theorem \ref{teorema_w(f)_totalmete_disconexo}, $\omega(f)$ is totally disconnected.
\end{proof1}

\begin{corollary}\label{E(X) closed}
Let $D$ be a dendrite with $\cl_D(E(D))=E(D)$. If $f:D\to
D$ is a transitive map, then $\omega(f)$ is totally disconnected.
\end{corollary}

\begin{proof1}
Take an arc $[a,b]$ contained in $D$, with $a\neq b$, we consider the following 2 cases.
\setcounter{case}{0}
\begin{case}
$(a,b)\cap R(D)=\emptyset$.
\end{case}
Then $[a,b]$ is a free arc in $D$.
 
\begin{case}
$(a,b)\cap R(D)\neq\emptyset$.
\end{case}
Take a point $x\in(a,b)\cap R(D)$. Since $D$ is a dendrite whose set of endpoints is closed, we know, by \cite[Corollary 3.6]{DWCE}, that $R(D)$ is discrete. Hence, we can take a connected open subset $V$ of $D$ such that $V\cap R(D)=\{x\}$. It follows that $(V\cap(a,b))\setminus\{x\}$  has exactly 2 components; the closure of either of them is a free arc in $D$. 

In either case, the arc $[a,b]$ contains a free arc. Therefore, by Corollary 
\ref{corolario_dendrita_totalmente_disconexo}, $\omega(f)$ is totally discconected.
\end{proof1}

\medskip
The following result is an immediate consequence of either Corollary \ref{corolario_dendrita_totalmente_disconexo} or Corollary \ref{E(X) closed}.

\begin{corollary}
Let $f:[0,1]\to[0,1]$ be a transitive map. Then $\omega(f)$ is totally disconnected.
\end{corollary}

\section{The example}
\label{example}
In several places of this section we consider sequences of the form $\lbrace n_{i} \rbrace \subset \mathbb{N}$. All of them are strictly increasing.
Let $\lbrace e_{n} \rbrace \subset \mathbb{R}$ be a sequence contained in the unit interval $[0,1]$. Define the \textit{omega limit set of the sequence $\{e_n\}$} as
\[
\omega(\lbrace e_{n} \rbrace)= \lbrace y\in [0,1] : \mbox{there exists} \ \lbrace  n_{i} \rbrace \subset \mathbb{N}, \ \lim _{i\rightarrow \infty} e_{n_{i}}=y \rbrace.
\]

Let $\Gamma = \lbrace a_{n} \rbrace = \lbrace a_{1},a_{2},a_{3},\ldots \rbrace = \left\lbrace \frac{1}{2},\frac{1}{2^{2}},\frac{3}{2^{2}}, \frac{1}{2^{3}},\frac{3}{2^{3}}, \frac{5}{2^{3}},\frac{7}{2^{3}}, \frac{1}{2^{4}}, \ldots \right\rbrace$ be the dyadic rational numbers contained in the open interval $(0,1)$.

Each $r\in (0,1]$ defines the following subsequence of $\lbrace a_{n} \rbrace$,
\[
[0,r]\cap \lbrace a_{n} \rbrace = \lbrace b_{1},b_{2},b_{3},\ldots \rbrace = \lbrace b_{k}=a_{n_{k}} \rbrace.
\]

\begin{remark}
The omega limit sets of those sequences are:
\begin{itemize}
\item[(i)] $\omega( \lbrace a_{n} \rbrace )= [0,1]$, and
\item[(ii)] for each $r\in (0,1]$, $\omega( \lbrace b_{k} \rbrace )= [0,r]$.
\end{itemize}
\end{remark}

From now on we recall the construction of the universal dendrite $D_{\infty}$, contained in the plane $\mathbb{R}^{2}$, and simultaneously we define a map $f:D_{\infty} \rightarrow D_{\infty}$ with the following properties:
\begin{itemize}
\item[(i)] $f$ is transitive in $D_{\infty}$.
\item[(ii)] The hyperspace $\omega(f)$ is not totally disconnected.
\end{itemize}

In the construction of the dendrite $D_{\infty}$ we follow some ideas given in \cite{Cont}, Section 10,37, page 181. In the description of the map $f:D_{\infty} \rightarrow D_{\infty}$ we follow some ideas from \cite{vero}, Section 3, page 55.

In the definition of $D_{\infty}$ we consider a sequence $\lbrace D_{n}\rbrace$ of dendrites, contained in $\mathbb{R}^{2}$, such that for each $n\in \mathbb{N}$, $D_{n}\subset D_{n+1}$. At the end, $D_{\infty} = \cl_{\mathbb{R}^2}(\cup_{n=1}^{\infty} D_{n})$.

\emph{Step 1}. Let $c\in \mathbb{R}^{2}$. By a \emph{star} with center $c$ and \emph{beams} $B_{j}$ we mean a countable union $S=\cup _{j=0}^{\infty}B_{j}$ of convex arcs $B_{j}$ in $\mathbb{R}^{2}$, each of them having $c$ as one end point, such that
\begin{itemize}
\item[(i)] $B_{j}\cap B_{k}=\lbrace c \rbrace$, whenever $j\neq k$, and
\item[(ii)] $\diam(B_{j})\rightarrow 0$ as $j\rightarrow \infty$.
\end{itemize}

For each $j\geq 0$, let $h_{j}: [0,1] \rightarrow B_{j}$ be a linear homeomorphism with $h_{j}(0)=c$. The sequence $\lbrace h_{j}( a_{n}): a_{n} \in \Gamma \rbrace$ is a copy of $\lbrace a_{n} \rbrace$ contained in $B_{j}$. We call it the sequence of dyadic points of $B_{j}$.

Let $D_{1}$ be a star such that $c$ is the origin of $\mathbb{R}^{2}$, $c= \textbf{0}$, $B_{0}=[-1,0]\times \lbrace 0 \rbrace$, $B_{1}=\lbrace 0 \rbrace \times [0,1]$, and for each $j\geq 2$, $B_{j}\subset \lbrace (x,y): x\geq 0, \ y\geq 0 \rbrace$.

For each $j\geq 2$, let $m_{j}$ be the slope of $B_{j}$. Let us assume that $0<m_{j+1}<m_{j}$, and $\lim\limits_{j\rightarrow \infty}m_{j}=0$.

In order to describe any point in $D_{1} \setminus \lbrace \textbf{0} \rbrace$ we 
consider a pair of numbers $(n,r)$, where $n\geq 0$ and $r\in (0,1]$. Let $x\in D_{1} \setminus \lbrace \textbf{0} \rbrace$, $x\sim (n,r)$ means $x\in B_{n}$ and $d(x,\textbf{0})= r \cdot\lo (B_{n})$. The origin is described by $\textbf{0}\sim (0)$.

The dendrite $D_{1}$ is strongly invariant under $f$. The restriction of $f$ to $D_{1}$ has the following properties:
\begin{itemize}
  \item[(i)] $f(B_{0})=\lbrace \textbf{0} \rbrace$. Thus, if $x\in B_{0}$, with $x\neq \textbf{0}$, and $y=f(x)$, then $x\sim (0,r)$ and $y\sim (0)$. In this case we represent the action of $f$ as: $(0,r)\rightarrow (0)$. Also, we have that $(0)\rightarrow (0)$.
  \item[(ii)] For every $n\in \mathbb{N}$, $f:B_{n} \rightarrow B_{n-1}$ is a linear homeomorphism. If $x\in B_{n}$, with $x\sim (n,r)$, and $y=f(x)$, then $y\sim (n-1,r)$. The action of $f$ is given by $(n,r)\rightarrow (n-1,r)$.
  \item[(iii)] For every point $x\in D_{1}$, $\omega(x,f)=\lbrace \textbf{0} \rbrace$.
  \end{itemize}
  
\emph{Step 2}. Let $n\geq 0$ and $a_{i}\in \Gamma$. For each point $x\in B_{n}$, with $x\sim (n,a_{i})$, form a star $S_{(n,a_{i})}$ with center in $x$ such that $S_{(n,a_{i})}\cap D_{1}=\lbrace x \rbrace$. The collection of all stars $\lbrace S_{(n,a_{i})} : n\geq 0, \ a_{i}\in \Gamma\rbrace$ satisfy the following properties:
\begin{itemize}
\item[(i)] For each $n,m \geq 0$, and each $a_{i},a_{j}\in \Gamma$, $S_{(n,a_{i})} \cap S_{(m,a_{j})}= \emptyset$ provided that $(n,a_{i}) \neq (m,a_{j})$.
\item[(ii)] If $n\rightarrow \infty$ or $i\rightarrow \infty$, then $\diam(S_{(n,a_{i})})\rightarrow 0$.
\end{itemize}
 
Let $D_{2}=D_{1}\cup(\cup \lbrace S_{(n,a_{i})} : n\geq 0, \ a_{i}\in \Gamma \rbrace)$. 

Let $x\in D_{2}\setminus D_{1}$. With $x\sim (n,a_{i},m,r)$ we mean that $x$ is a point in the $m$ beam of the star $S_{(n,a_{i})}$. We denote this with $x\in B_{m,(n,a_{i})}$. Let $p\in D_{1}$ be the point with $p\sim (n,a_{i})$. Then $d(x,p)= r \cdot \lo(B_{m,(n,a_{i})})$.

The dendrite $D_{2}$ is strongly invariant under $f$. The restriction of $f$ to dendrite $D_{2}$ has the following properties:
\begin{itemize}
  \item[(i)] Let $n\geq 1$.  
  $$f: B_{n}\cup (\cup \lbrace S_{(n,a_{i})} : a_{i}\in \Gamma \rbrace) \rightarrow B_{n-1}\cup (\cup \lbrace S_{(n-1,a_{i})} : a_{i}\in \Gamma \rbrace)$$ 
is a homeomorphism. Let $x\in D_{2}\setminus D_{1}$. If $x\sim (n,a_{i},m,r)$ and $y=f(x)$, then $y \sim (n-1,a_{i},m,r)$. In this case we represent the action of $f$ as: 
\[
(n,a_{i},m,r) \rightarrow (n-1,a_{i},m,r).
\]
  \item[(ii)] Let $x\in D_{2}\setminus D_{1}$. If $x\sim (0,a_{i},m,r)$ with $a_{i}=\frac{p}{2^{l}}$, consider the stars $S_{(0,a_{i})}$ and $\cup_{j=l}^{\infty}B_{j}$. The map $f:S_{(0,a_{i})} \rightarrow \cup_{j=l}^{\infty}B_{j}$ is a  homeomorphism. If $y=f(x)$, then $y\sim (m+l,r)$. In this case the action of $f$ is given by 
  $$(0,a_{i},m,r) \rightarrow (m+l,r) . $$
  Note that $f$ sends each star $S_{(0,\frac{k}{2^{l}})}$, with $1\leq k \leq 2^{l}-1$, $k$ odd number, to the star $\cup_{j=l}^{\infty}B_{j}$.
  \item[(iii)] For each point $x\in D_{2}$ there exists $N\in \mathbb{N}$ such that $f^{N}(x)= \textbf{0}$. Therefore, $\omega(x,f)=\lbrace \textbf{0} \rbrace$.
  \end{itemize}
  
\emph{Step 3}. Define the dendrite $D_{3}$ in a similar manner as $D_{2}$ was defined from $D_{1}$. At each point $x\in D_{2}\setminus D_{1}$, with $x\sim (n,a_{i},m,a_{j})$, with $n,m \geq 0$, $a_{i}, a_{j} \in \Gamma$, form a star $S_{(n,a_{i},m,a_{j})}$ with center in $x$, such that 
$S_{(n,a_{i},m,a_{j})}\cap D_{2}=\lbrace x \rbrace$.

The collection of all stars $\lbrace S_{(n,a_{i},m,a_{j})} : n,m \geq 0, \ a_{i}, a_{j} \in \Gamma\rbrace$ satisfy the following properties:
\begin{itemize}
\item[(i)] $S_{(n,a_{i},m,a_{j})} \cap S_{(s,a_{l},t,a_{k})}= \emptyset$, whenever $(n,a_{i},m,a_{j}) \neq (s,a_{l},t,a_{k})$.
\item[(ii)] If $n\rightarrow \infty$, $m\rightarrow \infty$, $i\rightarrow \infty$ or $j\rightarrow \infty$, then $\diam(S_{(n,a_{i},m,a_{j})})\rightarrow 0$.
\end{itemize}
 
Let $D_{3}=D_{2}\cup(\cup \lbrace S_{(n,a_{i},m,a_{j})} : n,m \geq 0, \ a_{i},a_{j} \in \Gamma \rbrace)$. 

The dendrite $D_{3}$ is strongly invariant under $f$. The restriction of $f$ to dendrite $D_{3}$ has the following properties:
\begin{itemize}
  \item[(i)] Let $n\geq 1$. Let $x\in D_{3}\setminus D_{2}$. If $x\sim (n,a_{i},m,a_{j},t,r)$ and $y=f(x)$, then $y \sim (n-1,a_{i},m,a_{j},t,r)$. We represent the action of $f$ as: 
\[
(n,a_{i},m,a_{j},t,r) \rightarrow (n-1,a_{i},m,a_{j},t,r).
\]
  \item[(ii)] Let $x\in D_{3}\setminus D_{2}$. If $x\sim (0,a_{i},m,a_{j},t,r)$ with $a_{i}=\frac{p}{2^{l}}$, and $y=f(x)$, then $y\sim (m+l,a_{j},t,r)$. In this case the action of $f$ is given by 
  $$(0,a_{i},m,a_{j},t,r) \rightarrow (m+l,a_{j},t,r) . $$
  \item[(iii)] For each point $x\in D_{3}$ there exists $N\in \mathbb{N}$ such that $f^{N}(x)= \textbf{0}$. Therefore, $\omega(x,f)=\lbrace \textbf{0} \rbrace$.
  \end{itemize}

\emph{Step 4}. For each $n\in \mathbb{N}$ define the dendrite $D_{n}$ following the directions given in steps 1, 2 and 3. Let $D_{\infty}$ be the closure of the union of the family of those dendrites, $D_{\infty}=cl(\cup_{n=1}^{\infty} D_{n})$. The continuum $D_{\infty}$ is the universal dendrite. 

Let $x\in (\cup_{n=1}^{\infty} D_{n})\setminus \lbrace \textbf{0} \rbrace$. Then
\[
x \sim (n_{1},a_{m_{1}},n_{2},a_{m_{2}},\ldots ,n_{k},r),
\]
for some $n_{j}\geq 0$, with $1\leq j \leq k$, and for some $a_{m_{j}}\in \Gamma$, with $1\leq j \leq k-1$, and some $0<r \leq 1$. We call the vector $(n_{1},a_{m_{1}},n_{2},a_{m_{2}},\ldots ,n_{k},r)$ the \emph{itinerary} of $x$.

\begin{remark}
\label{D-inifinita-props}
The following is a list of some properties of $D_{\infty}$, and of some properties of the map $f:D_{\infty} \rightarrow D_{\infty}$.

\begin{itemize}
\item[(i)] The set $D_{\infty}\setminus \lbrace \textbf{0} \rbrace$ has infinite countable components, $D_{\infty}\setminus \lbrace \textbf{0} \rbrace = \cup _{n=0}^{\infty}C_{n}$. For each $n\geq 0$ let $E_{n}=\cl_{D_\infty}(C_{n})$. For each $n\geq 0$, we have that $B_{n}\subset E_{n}$. For each pair $n\neq m$, $E_{n}\cap E_{m}=\lbrace \textbf{0} \rbrace$. And $\lim\limits_{n\rightarrow \infty} \diam(E_{n})=0$.
\item[(ii)] For each $n\geq 1$, $f:E_{n}\rightarrow E_{n-1}$ is a homeomorphism.
\item[(iii)] For each $k\in \mathbb{N}$, the union $\cup_{n=k}^{\infty}E_{n}$ is homeomorphic to $D_{\infty}$. 
\item[(iv)] $f$ sends $E_{0}$ onto $\cup_{n=1}^{\infty}E_{n}$.
\item[(v)] Let $i\in \mathbb{N}$. Let $y\in B_{0}$ be the point with itinerary $(0,a_{i})$. Let $E_{a_{i}}\subset D_{\infty}$ be the following dendrite
\[
E_{a_{i}}=\lbrace y\rbrace \cup \lbrace x\in D_{\infty}: x\neq y, \ [x,y]\cap B_{0}=\lbrace y \rbrace \rbrace.
\]
The dendrite $E_{a_{i}}$ is homeomorphic to $D_{\infty}$. 
\item[(vi)] Let $i\in \mathbb{N}$. Then $\diam(E_{a_{i}})\rightarrow 0$ as $i\rightarrow \infty$.
\end{itemize}
\end{remark}

\begin{proposition} 
\label{transitive}
The map $f:D_{\infty} \rightarrow D_{\infty}$ is transitive.
\end{proposition}

\begin{proof1}
In \cite{vero}, Theorem 3.7, it is proved that $f$ is mixing. That mixing implies transitivity is immediate from definitions.
\end{proof1}

\medskip

Let $x\in D_{\infty} \setminus (\cup_{n=1}^{\infty} D_{n})$. There exists a sequence $\lbrace x_{k} \rbrace \subset \cup_{n=1}^{\infty} D_{n}$ such that $\lim\limits _{k\rightarrow \infty}x_{k}=x$. For each $k \in \mathbb{N}$, the itineraries of $x_{k}$ and $x_{k+1}$ are:
\[
x_{k}\sim (n_{1},a_{m_{1}},n_{2},a_{m_{2}},\ldots , n_{k},a_{m_{k}}),
\]
and
\[
x_{k+1}\sim (n_{1},a_{m_{1}},n_{2},a_{m_{2}},\ldots , n_{k},a_{m_{k}}, n_{k+1},a_{m_{k+1}}).
\]

Define the itinerary of $x$ as this vector with countable many coordinates:
\[
x\sim (n_{1},a_{m_{1}},n_{2},a_{m_{2}},\ldots , n_{k},a_{m_{k}},n_{k+1},a_{m_{k+1}}, \ldots).
\]

Let $y=f(x)$. Then
\[
y\sim (n_{1}-1,a_{m_{1}},n_{2},a_{m_{2}},\ldots , n_{k},a_{m_{k}},n_{k+1},a_{m_{k+1}}, \ldots),
\]
if $n_{1}\geq 1$, and
\[
y\sim (n_{2}+l,a_{m_{2}},n_{3},a_{m_{3}},\ldots , n_{k},a_{m_{k}},n_{k+1},a_{m_{k+1}}, \ldots),
\]
if $n_{1}=0$, and $a_{m_{1}}=\frac{j}{2^{l}}$, with $1\leq j \leq 2^{l}-1$.

Let $j\geq 0$. For each value $0< r \leq 1$ consider the arc $[\textbf{0},u_{j}]\subset B_{j}$, where the distance from $u_{j}$ to $\textbf{0}$ is $r \cdot \lo(B_{j})$. Let $D_{r}=\cup_{j=0}^{\infty}[\textbf{0},u_{j}]$. Let $D_{0}=\lbrace \textbf{0} \rbrace$.

Let $\varphi : [0,1]\rightarrow 2^{D_{\infty}}$ be the function given by $\varphi(r)= D_{r}$, for each $r\in [0,1]$. It is not difficult to see that $\varphi$ is a homeomorphism from the unit interval $[0,1]$ to the collection $\lbrace D_{r}:r\in [0,1]\rbrace$. 

Hence, $\varphi([0,1])$ is an arc contained in the hyperspace $2^{D_{\infty}}$.

\begin{proposition}
Let $f:D_{\infty} \rightarrow D_{\infty}$ be the map described above. Then the hyperspace $\omega(f)$ is not totally disconnected.
\end{proposition}

\begin{proof1}
It is enough to show that for each $r\in (0,1]$, there exists a point $x\in D_{\infty}$ such that $\omega(x,f)=D_{r}$.

Let $r\in (0,1]$. Let $\lbrace b_{k} \rbrace$ be a subsequence of $\Gamma$ such that $\lbrace b_{k} \rbrace = \Gamma \cap [0,r]$. Here $b_{k}=a_{n_{k}}$, for each $k\geq 1$, and the sequence $\lbrace n_{k} \rbrace$ is strictly increasing.

For each $s\geq 0$, consider the arc $[\textbf{0},u_{s}]\subset B_{s}$, where  $u_{s}\sim (s,r)$.

Let $x\in D_{\infty}$ be the point with the following itinerary:
\[
x\sim (0,b_{1},0,b_{2},0,b_{3},\ldots).
\]

\textbf{Claim.} $\omega(x,f)=D_{r}$.

Let $j,l \in \mathbb{N}$ such that $b_{1}=\frac{j}{2^{l}}$. Note the following: the itineraries of the first points in the orbit $o(x,f)$ are:
\[
(0,b_{1},0,b_{2},0,b_{3},\ldots) \ \rightarrow \ (l,b_{2},0,b_{3},0,b_{4},\ldots)
\]
\[
\rightarrow \ (l-1,b_{2},0,b_{3},0,b_{4},\ldots) \ \rightarrow \ (l-2,b_{2},0,b_{3},0,b_{4},\ldots)
\]
\[
\rightarrow \cdots \rightarrow \ (0,b_{2},0,b_{3},0,b_{4},\ldots)
\]

For each $k\geq 2$, there exists $m_{k}\in \mathbb{N}$ such that $y=f^{m_{k}}(x)$ with
\[
y \sim (0,b_{k},0,b_{k+1},0,b_{k+2},\ldots).
\]

For each $k\in \mathbb{N}$, let $e_{k}\in B_{0}$ be the point that $e_{k}\sim (0,b_{k})$. According to the claim (vi) of Remark \ref{D-inifinita-props}, the distance $d(f^{m_{k}}(x),e_{k}) \rightarrow 0$ as $k\rightarrow \infty$. It follows that the arc $[\textbf{0},u_{0}]$ is contained in the limit set $\omega(x,f)$. 

Now, the set $\omega(x,f)$ is strongly invariant under map $f$, hence for every $s\in \mathbb{N}$, $[\textbf{0},u_{s}]\subset \omega(x,f)$. Thus, $D_{r}\subset \omega(x,f)$.

Let $x_{0}\in E_{0} \setminus [\textbf{0},u_{0}]$. Let $\delta_{1}=\dist(x_{0}, \cup_{n=1}^{\infty}E_{n})$, and let $\delta_{2}=\dist(x_{0}, [\textbf{0},u_{0}])$.
Let $z\in o(x,f)$. If the first coordinate of the itinerary of $z$ is distinct from 0, then $z\in \cup_{n=1}^{\infty}E_{n}$, hence $d(x_{0},z)\geq \delta_{1}$. If the first coordinate of $z$ is 0, there exists $k\in \mathbb{N}$ such that $z=f^{m_{k}}(x)$. Since $d(f^{m_{k}}(x),e_{k}) \rightarrow 0$ as $k\rightarrow \infty$, there exists $K\in \mathbb{N}$ such that for each $k\geq K$, $d(x_{0},f^{m_{k}}(x))\geq \frac{\delta_{2}}{2}$. It follows that $x_{0}\notin \omega(x,f)$.
Using again that $\omega(x,f)$ is a strongly invariant set under map $f$, we conclude that $\omega(x,f)=D_{r}$.
\end{proof1}

\end{document}